\documentclass[12pt]{amsart}

\newcommand {\C}{\mbox{l\hspace{-.47em}C}} 
\newcommand {\D}{\mbox{I\hspace{-.15em}D}} 

\newtheorem{thm}{Theorem}[section]

\newtheorem{lem}[thm]{Lemma}
\newtheorem{cor}[thm]{Corollary}

\begin{document}

\title[Compressions of Resolvents]{Compressions of Resolvents 
and Maximal \\
Radius of Regularity}
\author{C. BADEA}
\address{URA 751 au CNRS \&
  UFR de Math\'ematiques, Universit\'e de Lille I,
  F--59655 Villeneuve d'Ascq, France}
\email{badea@@gat.univ-lille1.fr}
\author{
		M. MBEKHTA}
\address{URA 751 au CNRS \&
  UFR de Math\'ematiques, Universit\'e de Lille I,
  F--59655 Villeneuve d'Ascq, France}
\curraddr{University of Galatasaray, \c{C}iragan Cad no 102,
		Ortakoy 80840, Istanbul, Turkey}
\email{mbekhta@@gat.univ-lille1.fr}


\begin{abstract}
Suppose that $\lambda - T$ is left-invertible in $L(H)$ for all
$\lambda \in \Omega$, where $\Omega$ is an open subset of the 
complex plane. Then an operator-valued function $L(\lambda)$ is a left 
resolvent of $T$ in $\Omega$ if and only if $T$ has an extension 
$\tilde{T}$,
the resolvent of which is a dilation of $L(\lambda)$ of a particular form. 
Generalized resolvents exist on every open set $U$, with $\overline{U}$ 
included 
in the regular domain of $T$. This implies a formula for the maximal 
radius of regularity of $T$ in terms of the spectral radius of its 
generalized inverses. A solution to an open problem raised 
by J. Zem\'anek is obtained. 
\end{abstract}

\maketitle


\section{Introduction}
Let $H$ be a complex, separable, infinite dimensional  Hilbert space and 
let 
$L(H)$ be the C$^{\ast}$-algebra of all 
continuous
linear operators on $H$. The resolvent set $\rho(T)$ of $T \in L(H)$ 
is, by definition, 
the set of all complex numbers $\lambda \in \C$ such that the operator 
$\lambda I - T$ is invertible in the algebra $L(H)$. Then its resolvent 
$R(\lambda) = (\lambda I - T)^{-1}$ is an analytic function on $\rho(T)$ 
such 
that $R(\lambda)(\lambda I - T) = (\lambda I - T)R(\lambda) = I$ for all 
$\lambda \in \rho(T)$. Here 
$I$ is the identity operator on $H$. 
Moreover, $R(\lambda)$ satisfies the first 
resolvent equation, namely
$$R(\lambda) - R(\mu) = (\mu - \lambda)R(\lambda)R(\mu)$$
for all $\lambda , \mu \in \rho(T)$.
The spectrum $\sigma(T)$ of 
$T$ is the complement of $\rho(T)$ in $\C$. 

The left resolvent set $\rho_{\ell}(T)$ of the bounded linear operator $T$ 
is 
defined as the set of all complex numbers $\lambda \in \C$ such that 
$\lambda I - T$ is left invertible in $L(H)$, that is one-to-one and with 
closed range. The 
left spectrum of $T$ is the set $\sigma_{\ell}(T) = \C \setminus 
\rho_{\ell}(T)$. 
According to a result due to 
G. R. Allan \cite{All1, All2}, there exists an analytic function 
$L(\lambda)$ on $\rho_{\ell}(T)$ such that $L(\lambda)(\lambda I - T) = I$,
for all $\lambda \in \rho_{\ell}(T)$. An operator-valued mapping
$L(\lambda) : U \to L(H)$, $U \subseteq \rho_{\ell}(T)$, is said to be 
a {\sl left resolvent} of $T$ on $U$,  
if $L(\lambda)(\lambda I - T) = I$ on $U$ and, in addition, 
$L(\lambda)$ satisfies the resolvent equation above for all $\lambda$ and 
$\mu$ in the {\sl same} connected component of $U$. If $U = 
\rho_{\ell}(T)$, 
the above $L(\lambda)$ is simply called a {\sl left resolvent} of $T$. 
The right 
resolvent set $\rho_{r}(T)$ and right resolvents are defined in a 
similar way. 

Let $reg(T)$ denote the set of all complex numbers $\lambda$ for which 
$\mu I - T$ possess an {\sl{analytic generalized inverse}}, for all $\mu$ 
in a 
neighborhood of $\lambda$, that is there exists $U_{\lambda}$ a 
neighborhood of 
$\lambda$ and an operator-valued analytic function 
$G(\mu) : U_{\lambda} \to L(H)$ such that, for all $\mu \in U_{\lambda}$,
\begin{equation}
(\mu I - T)G(\mu)(\mu I - T) = (\mu I - T)
\end{equation}
and
\begin{equation}
G(\mu)(\mu I - T)G(\mu) = G(\mu) .
\end{equation}
Then $\sigma_{g}(T) = \C \setminus reg(T)$ is the {\sl{generalized 
spectrum}} of $T$ 
\cite{Mbe1}, \cite{Mbe2}. Note also that $reg(T)$ has several equivalent 
definitions. 
According to a result of Shubin \cite{Shu}, there 
exists a global analytical function $G(\lambda)$ defined on $reg(T)$ which 
is 
a generalized inverse of $\lambda I - T$ for all $\lambda \in reg(T)$. By 
definition, a 
{\sl{generalized resolvent}} of $T$ on $U$, $U \subseteq reg(T)$, 
is an operator-valued mapping
$G(\lambda)$ satisfying the conditions (1) and (2) for $\mu \in U$ and the
resolvent equation for all $\lambda$ and $\mu$ in the same connected 
component of $U$. {\sl Generalized resolvents} of $T$ 
are generalized resolvents of $T$ on $reg(T)$. 
If $U \subseteq \rho (T)$, then a generalized (left, right) resolvent on 
$U$ is 
unique and coincides with the usual resolvent.  

It is not known if left or generalized resolvents always exist. These open 
problems go back to C. Apostol and K. Clancey \cite{ApCl1}, \cite{ApCl2}
and P. Saphar \cite{Sap1}, \cite{Sap2}. See also the list of references in 
\cite{BaMb} for several contributions on these problems.

The aim of the present paper is twofold. Firstly, we will prove that 
generalized resolvents 
for $T$ 
always exist on every open subset $U$ of $reg(T)$ such that $\overline{U} 
\subset reg(T)$. This implies a formula for the maximal radius 
of regularity of $T$ in terms of the spectral radius of its generalized 
inverses. 
In particular, one obtains (section 2), in the case of Hilbert 
space operators, a positive answer to a problem recently 
raised by Zem\'anek \cite{Zem} for Banach 
algebra elements. Note that in \cite{BaMb} several partial answers to 
Zem\'anek's question (for closed operators with a dense domain on a 
Banach space) were given.

Secondly, we give a new characterization of left resolvents 
in terms of dilations and compressions. The result of C.J. Read 
\cite{Rea3},
about the extent to which an extension of a Hilbert space operator reduces 
the spectrum, implies that every operator $T$ admits 
an extension $\tilde{T}$ to a larger Hilbert space $\tilde{H} = H \oplus 
H'$ 
such that 
$\lambda I - \tilde{T}$ is invertible in $L(\tilde{H})$ for all 
$\lambda \in \rho_{\ell}(T)$. Moreover, an analytic left inverse function 
$L(\lambda)$ of $\lambda I - T$ can be obtained as the compression on $H$ 
of $R(\tilde{T})(\lambda) = (\lambda I - \tilde{T})^{-1}$. 
The matrix of $R(\tilde{T})(\lambda)$ with respect to the 
decomposition $\tilde{H} = H \oplus H'$ will then be of the form

$$\left(\begin{array}{cc}
			\displaystyle{L(\lambda)} & \ast \\
			\ast & \ast
		\end{array}\right) \ ,$$

\noindent where $\ast$ means a suitable entry.

The main result of the second part of this note is to prove that 
$L(\lambda)$ is a left resolvent
of $T$ in $\Omega \subseteq \rho_{\ell}(T) \cap \sigma(T)$ if and only if 
there exists an 
extension $\tilde{T}$ of $T$ on a larger Hilbert space 
$\tilde{H} = H \oplus H'$ such that 
$ R(\tilde{T})(\lambda) = (\lambda I - \tilde{T})^{-1}$ is invertible in 
$L(\tilde{H})$ for all 
$\lambda \in \Omega$ and the matrix of $R(\tilde{T})(\lambda)$ with 
respect to the 
decomposition $\tilde{H} = H \oplus H'$ has the form   

$$\left(\begin{array}{cc}
			\displaystyle{L(\lambda)} & 0 \\
			\ast & \ast
		\end{array}\right) \ ,$$

\noindent that is $R(\tilde{T})(\lambda)$ leaves $H'$ invariant, 
for all $\lambda \in \Omega$. In fact, a stronger version of 
this result will be proved.  

\section{Generalized resolvents on subsets of $reg(T)$ and radius of 
regularity}

\subsection{Generalized resolvents on subsets of $reg(T)$}

The following result is a partial result for the one-sided resolvent 
problem.

\begin{lem}[\cite{HerII}]
Let $T\in L(H)$ and let $\Omega \subset \overline{\Omega} \subset 
\rho_{\ell}(T)$ be an 
open set. Then there is a left resolvent of $T$ on $\Omega$.
\end{lem}

\noindent {\sl Proof. } The result follows from Proposition 9.17 
\cite{HerII} ; see 
also Lemma 3.18 \cite{HerI} and Errata to \cite{HerI}.
\qed
\bigskip

The following result is the corresponding extension for generalized 
resolvents.

\begin{thm}
Let $T\in L(H)$ and let $\Omega \subset \overline{\Omega} \subset reg(T)$ 
be an 
open set. Then there is a generalized resolvent of $T$ on $\Omega$, 
that is there exists an analytic function 
$G(\lambda) : \Omega \to L(H)$ such that $G(\lambda)$ satisfies the 
resolvent 
identity for all $\lambda , \mu$ in the same connected component of $U$ 
and 
$$(\lambda I - T)G(\lambda)(\lambda I - T) = (\lambda I - T) , \, \, \, 
\lambda \in \Omega$$ 
as well as
$$G(\lambda)(\lambda I - T)G(\lambda) = G(\lambda) , \, \, \, \lambda \in 
\Omega .$$
\end{thm}

\noindent {\sl Proof. } Let $H = H_r \oplus H_0 \oplus H_l$ be the Apostol 
decomposition 
\cite{Apo}, \cite{LaMb}. With respect to this decomposition, the operator 
$T$ 
can be written as 
$$\left(\begin{array}{ccc}
			\displaystyle{T_r} & A & B \\
			0 & T_0 & C \\
			0 & 0 & T_l
		\end{array}\right) \  ,$$
for suitable operator entries. Then \cite[Th\'eor\`eme 4.10]{LaMb} 
$$reg(T) = \rho_r(T_r) \cap \rho(T_0) \cap \rho_{\ell}(T_l).$$
Using twice Lemma 2.1, we find a left resolvent $L(\lambda)$ for $T_l$ and 
a right 
resolvent $R(\lambda)$ for $T_r$, both defined on 
$\Omega \subset \overline{\Omega} \subset reg(T)$. In particular, 
$(\lambda I - T_r)R(\lambda) = I$ and $L(\lambda)(\lambda I - T_{\ell}) = 
I$ for all $\lambda \in \Omega$. 
Consider $G : \Omega \to 
L(H)$, given by the following matrix, with respect to the Apostol 
decomposition :
$$\left(\begin{array}{ccc}
			\displaystyle{R(\lambda)} & -R(\lambda)AR(T_0)(\lambda) & 
R(\lambda)[AR(T_0)(\lambda)C - B]L(\lambda) \\
			0 & R(T_0)(\lambda) & -R(T_0)(\lambda)CL(\lambda) \\
			0 & 0 & L(\lambda)
		\end{array}\right) \  ,$$
where $R(T_0)(\lambda) = (\lambda I - T_0)^{-1}$. Denote $P(\lambda) = 
(\lambda I - T)G(\lambda)$ and 
$Q(\lambda) = G(\lambda)(\lambda I - T)$ for $\lambda \in \Omega$. Then 
simple computations show that 
$$P(\lambda) = \left(\begin{array}{ccc}
			I &0 & 0 \\
			0 &I & 0 \\
			0 & 0 & \displaystyle{(\lambda I - T_{\ell})L(\lambda)}
		\end{array}\right) \  $$
and
$$Q(\lambda) = \left(\begin{array}{ccc}
			\displaystyle{R(\lambda)(\lambda I - T_{r})} &0 & 0 \\
			0 &I & 0 \\
			0 & 0 & I
		\end{array}\right) \  .$$
 
Therefore $G(\lambda)$ is a global generalized 
inverse for $T$ on $\Omega$. We show now that $G(\lambda)$ 
satisfies the resolvent identity on $\Omega$. To this end, we use a 
criterion given in 
\cite[Theorem 2.7]{BaMb}. It is easy to show that the range of 
$P(\lambda)$ coincides with 
the range of $\lambda I - T$ and the kernel of $Q(\lambda)$ coincides with 
the kernel of $\lambda I - T$. 
On the other hand, one has 
$$P(\lambda)P(\mu) =  
\left(\begin{array}{ccc}
			I &0 & 0 \\
			0 &I & 0 \\
			0 & 0 & \displaystyle{(\lambda I - T_{\ell})L(\lambda)(\mu I - 
T_{\ell})L(\mu)}
		\end{array}\right) \   .$$
But, using the fact that $L(\lambda)$ is a left resolvent of $\lambda I - 
T_{\ell}(T)$, we obtain 
\begin{eqnarray*}
(\lambda I - T_{\ell})L(\lambda)(\mu I - T_{\ell})L(\mu) & = & 
(\lambda I - T_{\ell})L(\lambda)(\mu I - \lambda I + \lambda I - 
T_{\ell})L(\mu) \\
     & = &  
(\lambda I - T_{\ell})(\mu - \lambda)L(\lambda)L(\mu) + (\lambda I - 
T_{\ell})L(\mu) \\
  & = &  
(\lambda I - T_{\ell})(L(\lambda) - L(\mu)) + (\lambda I - T_{\ell})L(\mu) 
\\
  & = & (\lambda I - T_{\ell})L(\lambda).
\end{eqnarray*}
Therefore $P(\lambda)P(\mu) = P(\lambda)$ and, similarly, 
$Q(\lambda)Q(\mu) = Q(\mu)$. 
Using \cite[Proof of Theorem 2.7]{BaMb}, we get that $G(\lambda)$ 
satisfies the resolvent 
identity in $\Omega$.

\qed
\bigskip

\subsection{Maximal radius of regularity}

As an application of the existence of generalized resolvents, the 
following formula 
for the maximal radius of regularity can be proved. It is a generalization 
of \cite[Theorem 3.1]{BaMb}.

\begin{thm}\label{genZem}
Let $T \in L(H)$ be a linear operator such that $0 \in reg(T)$.
Then 
$${\rm{ dist}} (0,\sigma_{g}(T)) = \sup \{ \frac{1}{r_{\sigma}(S)} : TST = 
T\},$$
where $r_{\sigma}(S)$ is the spectral radius of $S$.
\end{thm}  

{\sl{Proof. }} The proof is similar to that of \cite[Theorem 3.1]{BaMb}, 
where 
the present theorem has been proved (for closed, densely defined, Banach 
space 
operators) under the additional assumption that 
$T$ is Fredholm.  Therefore some details will be omitted. 
Let $S \in L(H)$ be a generalized inverse of
$T$. Then, using for instance \cite[Theorem 2.4]{BaMb}, we get 
$T^{n}S^{n}T^{n} = T^{n}$, for all
$n \geq1$. Let $\gamma (T)$ be the reduced minimum modulus
of $T$ : 
$$\gamma (T) = \inf \{ \| Tx \| : \mbox{ dist} (x,N(T)) = 1\},$$
where $N(T)$ is the kernel of $T$. By \cite[Lemme 3.5]{Mbe2}, 
we have $\gamma (T^{n}) \geq 1/(\| S^{n} \|)$ and
therefore
$$\lim_{n \to\infty} \gamma (T^{n})^{1/n} \geq \frac{1}{r_{\sigma}(S)}.$$
Using \cite[Th\'eor\`eme 3.1]{Mbe2} we get the inequality
$${\rm{ dist}}(0,\sigma_{g}(T)) \geq \sup\{ \frac{1}{r_{\sigma (S)}} : TST 
= T\} .$$

In order to prove the other inequality, set $d = {\rm{ 
dist}}(0,\sigma_{g}(T))$.
Then
$$B(0,d) = \{ \lambda \in {\bf C} : |\lambda | < d\} \subseteq reg(T).$$

Let $\varepsilon$ be a positive number and put
$$U = B(0,\frac{d}{1+\varepsilon}) = \{\lambda \in \C : 
|\lambda | < \frac{d}{1+\varepsilon} \} .$$
Then $U \subset \overline{U} \subset reg(T)$. Using Theorem 2.2,
there is a generalized resolvent $G(\lambda)$ for $T$ on $U$. Then,
for all $\lambda \in U$,
$$G(\lambda) = \sum_{n=0}^{\infty}\lambda ^{n}T_{n} , \quad T_{n} \in L(H) 
.$$
As in \cite[Proof of Theorem 3.1]{BaMb}, we can use the Cauchy's integral 
formula and 
\cite[Theorem 2.6]{BaMb} to get 
$$\| T_{n}\| \leq M\left( \frac{1+\varepsilon}{d}\right)^{n+1}, n\geq 0 ,$$
where $M = max\{ \| G(\lambda)\| : \lambda \in \overline{U}\}$, and 
$T_n = (-1)^{n+1}T_0^{n+1}$. 
Thus
$$\| T_{0}^{n+1}\| \leq M\left( \frac{1+\varepsilon}{d}\right)^{n+1}, 
n\geq 0,$$
which implies $r_{\sigma}(T_{0}) \leq (1+\varepsilon)/d$. Since $TT_{0}T = 
T$,
we have
$$\sup \{\frac{1}{r_{\sigma}(S)} : TST = T\} \geq \frac{1}{r_{\sigma} 
(T_{0})}
\geq \frac{d}{1+\varepsilon}.$$
Since $\varepsilon > 0$ was arbitrarily chosen, we have
$$\sup \{\frac{1}{r_{\sigma}(S)} : TST = T\} \geq d .$$
The proof is now complete. 
\qed

\bigskip

\noindent {\bf Remark.} The above proof shows that the $\sup$ is attained 
(for some $S_0$) 
in the 
formula 
$$ d := \mbox{ dist}(0,\sigma_{g}(T)) = \sup \{ \frac{1}{r_{\sigma}(S)} : 
TST = T\},$$
if one is able to construct a generalized resolvent for $T$ on $B(0,d)$. 
Conversely, if 
there exists $S_0 \in L(H)$ such that $TS_0 T = T$ and $d = 
\frac{1}{r_{\sigma}(S_0)}$, then 
$G(\lambda) = \sum_{n=0}^{\infty}\lambda ^{n}S_{0}^{n+1}$ is a generalized 
resolvent for $T$ on $B(0,d)$. We omit here the details.
\qed
\bigskip

Set $s(T) = \sup \{\frac{1}{r_{\sigma}(S)} : TST = T\}$

\begin{cor}
Let $T \in L(H)$ and suppose that $0 \in reg(T)$. Then we have $s(T^{n}) = 
s(T)^{n}$, 
for all
$n \geq 1$.
\end{cor}

{\sl{Proof }} Let $n \geq 1$.
Using the preceding Theorem and \cite[Theorem 3.1]{Mbe2}, 
we get
$$s(T) = \lim_{k\to \infty} \gamma (T^{k})^{1/k}.$$
Therefore $$s(T^{n}) = \lim_{k\to \infty} (\gamma (T^{kn})^{1/kn})^{n} 
= s(T)^{n}.$$ 
The proof is complete. 
\qed
\bigskip

The following is a solution for a problem raised by J. Zem\'anek 
\cite{Zem} in 
the more general setting of Banach algebras.

\begin{cor}
Let $T \in L(H)$ be a linear operator such that $0 \in \rho_{\ell}(T)$.
Then 
$${\rm{ dist}}(0,\sigma_{\ell}(T)) = \sup \{ \frac{1}{r_{\sigma}(S)} : ST 
= I\}.$$
\end{cor}

{\sl{Proof }}. The result follows from Theorem 2.3 and is similar to 
the proof of Corollary 3.4 in \cite{BaMb}.
\qed

\bigskip

Corollary 2.5 can be viewed as a one-sided generalization of the 
known formula 
$$0 \in \rho(T) \quad \Longrightarrow \quad {\rm{ dist}} (0, \sigma (T) ) 
= 
\frac{1}{r_{\sigma}(T^{-1})} .$$
\medskip

\section{Left resolvents as compressions of resolvents}

\subsection{Read's extension theorem}
The following result is a consequence of a result of Read \cite{Rea3} and 
Corollary 2.5.

\begin{lem}
Let $T \in L(H)$ and $\Omega \subseteq \rho_{\ell}(T)$. There 
exists a larger Hilbert space $\tilde{H} = H \oplus H'$ containing $H$  
and an 
extension $\tilde{T}$ of $T$ on $\tilde{H}$ such that $\lambda I - 
\tilde{T}$
is invertible in $L(\tilde{H})$ for all $\lambda \in \Omega$. Moreover,
$$r_{\sigma}((\lambda I - \tilde{T})^{-1}) = \inf \{r_{\sigma}(S) : S(T - 
\lambda I) = I\}$$
for all $\lambda \in \Omega$.
\end{lem}

{\sl Proof}. According to a result due to Read \cite{Rea3}, there exist 
a larger Hilbert space $\tilde{H} = H \oplus H'$ and an extension 
$\tilde{T}$ 
of $T$ such that the spectrum $\sigma(\tilde{T})$ of $\tilde{T}$ is equal 
to the approximate point spectrum $\sigma_{ap}(T)$ 
of $T$, that is
$$\sigma(\tilde{T}) = \sigma_{ap}(T) := \{ \lambda \in {\bf C} : 
\inf \{ \| (\lambda I - T)x \| : x \in H , \| x \| = 1\} = 0 \} \ .$$
This coincides with the left spectrum $\sigma_{\ell}(T)$ for Hilbert space 
operators. 

For the second part, using Corollary 2.5, we have 
\begin{eqnarray*}
\frac{1}{r_{\sigma}((\lambda I - \tilde{T})^{-1})} & = & 
 dist(\lambda,\sigma(\tilde{T})) \\
  &  = & dist(\lambda,\sigma_{\ell}(T)) \\
  & = & \sup \{ \frac{1}{r_{\sigma}(S)} : S(T - \lambda I) = I\}.
\end{eqnarray*}
This implies the above equality.
\qed     
\bigskip

\noindent {\bf Remark } In Read's theorem \cite{Rea3}, $H'$ can be chosen 
as a copy of 
$H$ and the norm of $\tilde{T}$ almost the norm of $T$, that is, if 
$\varepsilon > 0$, 
$\tilde{T}$ (and $H' \cong H$) can be chosen such that $\| \tilde{T} \| 
\leq (1 + \varepsilon ) \| T \|$. 
We also want to note that $H' \cong H$ depends upon $T$ in a sensitive 
way. Indeed, 
consider an analytic, operator-valued function $ \Omega \ni \lambda \to 
T(\lambda) \in L(H)$. 
By Read's proof, it is possible to find a copy $H'(\lambda)$ of $H$ and 
extensions 
$\Omega \ni \lambda \to \tilde{T}(\lambda) 
\in L(H\oplus H'(\lambda))$ such that $\sigma (\tilde{T}(\lambda)) = 
\sigma_{\ell}(T(\lambda))$. 
However, there are examples where one cannot finds a {\sl universal} 
copy of $H$ in place of $H'(\lambda)$ and {\sl analytic} familly 
$\Omega \ni \lambda \to \tilde{T}(\lambda) 
\in L(H\oplus H')$ of extensions. We use for this purpose an example of 
Ransford \cite{Ran}. 
Consider $\Omega = \C , H = \ell^{2}$ and 
$$T(\lambda)(a_1, a_2, \dots ) = (a_1, \lambda a_1, a_2, \lambda a_2, 
\dots )$$ 
which is holomorphic. Then it can be shown that 
$$\sigma_{\ell}(T(\lambda)) 
\subseteq \{z : |z| = (1 + |\lambda |^2)^{1/2} \} .$$
Suppose that an analytic family $\tilde{T}(\lambda) \in L(H\oplus H)$, 
extensions of $T(\lambda)$, 
such that 
$\sigma (\tilde{T}(\lambda)) = \sigma_{\ell}(T(\lambda))$ would exists. 
Then $\lambda \to 
\sigma_{\ell}(T(\lambda))$ would be an analytic set-valued function as 
the spectrum of an analytic family of operators (cf. 
\cite{Slo} for definitions and properties of set-valued analytic 
functions). 
This is in contradiction with the fact that 
$$\phi(\lambda) = \sup \{ - \log |z| : z \in \C \setminus 
\rho_{\ell}(T(\lambda)) \} = 
- \log (1 + |\lambda |^2)^{1/2}$$
is not subharmonic on $\C$ (cf. \cite{Ran}). Indeed, $\phi$ attains a 
maximum at 0.

We wish to thank Dan Timotin for useful discussions concerning Read's 
theorem.

\subsection{Compressions of resolvents}

Lemma 3.1 implies Allan's result mentioned in Introduction. Indeed, let 

$$L(\lambda)h = P_{H}R(\tilde{T})(\lambda)h \ , h \in H$$

\noindent be the compression on $H$ of 
$R(\tilde{T})(\lambda) = (\lambda I - \tilde{T})^{-1}$. Then, using the 
equality 
$R(\tilde{T})(\lambda)(\lambda I - \tilde{T}) = I$ and the fact that 
$\tilde{T}$ is 
an extension on $\tilde{H}$ of $T$, we obtain
$$L(\lambda) (\lambda I - T) = I \  \mbox{on} \ H \ .$$ 
Since $L(\lambda)$ is analytic, we obtain a global, analytic 
left inverse function of $\lambda I - T$, i.~e. Allan's result. The matrix 
of $R(\tilde{T})(\lambda)$ with respect to the 
decomposition $\tilde{H} = H \oplus H'$ will then be of the form

$$\left(\begin{array}{cc}
			\displaystyle{L(\lambda)} & \ast \\
			\ast & \ast
		\end{array}\right) \ .$$

The following characterization of left resolvents is the main result of 
this section. 

\begin{thm}\label{dil}
Let $T \in L(H)$ and let $\Omega$ be an open, connected 
subset of $\rho_{\ell}(T) \cap \sigma(T)$. Then 
$L(\lambda)$ is a left resolvent
of $T$ in $\Omega$ if and only if there exists an 
extension $\tilde{T}$ of $T$ on a larger Hilbert space 
$\tilde{H} = H \oplus H'$ such that 
$ R(\tilde{T})(\lambda) = (\lambda I - \tilde{T})^{-1}$ exists in 
$L(\tilde{H})$ for all 
$\lambda \in \Omega$ and the matrix of $R(\tilde{T})(\lambda)$ with 
respect to the 
decomposition $\tilde{H} = H \oplus H'$ has the form   

$$\left(\begin{array}{cc}
			\displaystyle{L(\lambda)} & 0 \\
			T(\lambda) & V(\lambda)
		\end{array}\right) \ ,$$
with suitable operator-valued functions $T(\lambda)$ and $V(\lambda)$.
\end{thm} 

The condition $\Omega\subset \sigma(T)$ is explained by the fact that on 
$\Omega \cap 
\rho(T)$ the usual resolvent is also a left resolvent. The condition of 
connectedness of $\Omega$ is justified by the definition of left 
resolvents as 
left inverses satisfying the resolvent identity on every connected 
component.

For the ''only if`` part of Theorem \ref{dil}, the following 
stronger result can be proved.

\begin{thm}\label{dil+}
Let $T \in L(H)$ and let $\Omega \subset \rho_{\ell}(T)  \cap  \sigma(T)$ 
be an 
open, connected subset. Suppose 
that there exists an 
extension $\tilde{T}$ of $T$ on a larger Hilbert space 
$\tilde{H} = H \oplus H'$ such that 
$\lambda I - \tilde{T}$ is invertible in $L(\tilde{H})$ for all 
$\lambda \in \Omega$ and the matrix of $R(\tilde{T})(\lambda) = (\lambda I 
- \tilde{T})^{-1}$ 
with respect to the 
decomposition $\tilde{H} = H \oplus H'$ has the form   

$$\left(\begin{array}{cc}
			\displaystyle{L(\lambda)} & S(\lambda) \\
			T(\lambda) & V(\lambda)
		\end{array}\right) \ ,$$
for some operator-valued functions $L(\lambda), (S(\lambda)$ defined on 
$H$ and $H'$,respectively, with values in H and $T(\lambda), V(\lambda)$ 
defined on $H$ and $H'$, respectively, with values in $H'$.
Suppose that $S(\lambda)T(\lambda) = 0$ for all $\lambda \in \Omega$. Then 
$L(\lambda)$ is a left resolvent
of $T$ in $\Omega$.
\end{thm}

Clearly, the above condition is satisfied if $S(\lambda) \equiv 0$.

For the ''if`` part of Theorem \ref{dil}, we need a construction 
of independent interest. 
Namely, we will prove the following result. 

\begin{thm} \label{ext}
Suppose that $T \in L(H)$ and let $\Omega$ be an open, connected 
subset of $\rho_{\ell}(T) \cap \sigma(T)$. Let $L(\lambda)$ be a global, 
analytic left-inverse 
function of $T - \lambda I$. There exists a larger Hilbert space 
$\tilde{H} = H \oplus H'$ and a 
family of operators $\tilde{T}(\lambda) \in L(\tilde{H}), \lambda \in 
\Omega$, such that 
\begin{enumerate}
\item $\tilde{T}(\lambda)$ is an extension of $T$, for all $\lambda \in 
\Omega \;$ ;
\item $\lambda I - \tilde{T}(\lambda)$ is invertible, for all $\lambda 
\in  \Omega \;$ ;
\item With respect to the decomposition $\tilde{H} = H \oplus H'$, we have
$$ (\lambda I - \tilde{T}(\lambda))^{-1} = \left(\begin{array}{cc}
			\displaystyle{L(\lambda)} & 0 \\
			\ast & \ast
		\end{array}\right) \ $$
for suitable entries $\ast \;$ ;
\item $\tilde{T}(\lambda) \equiv \tilde{T}$ is independent of $\lambda$ if 
and only if 
$L(\lambda)$ is a left resolvent for $T$ on $\Omega \;$ ;
\item If $\dim N(L(\lambda)) < \infty$, then the map $\lambda \to 
\tilde{T}(\lambda)$ is 
analytic.
\end{enumerate}
\end{thm}

Theorem \ref{dil} follows from Theorem \ref{dil+} and Theorem \ref{ext}. 

Using Lemma 2.1 and Theorem \ref{dil}, we obtain the following 
consequence, 
improving part of Lemma 3.1.

\begin{cor}
Let $T \in L(H)$ and let $\Omega \subset \overline{\Omega} \subset 
\rho_{\ell}(T) 
\cap \sigma(T)$ be an open, connected set. Then there exists an extension 
$\tilde{T}$ of 
$T$ on a larger Hilbert space $\tilde{H} = H \oplus H'$ such that 
$R(\tilde{T})(\lambda) 
= (\lambda I - \tilde{T})^{-1}$ exists in $L(\tilde{H})$ for all $\lambda 
\in \Omega$ and 
$R(\tilde{T})(\lambda)(H') \subseteq H'$, for all $\lambda \in \Omega$.
\end{cor} 


\subsection{Proof of Theorem \ref{dil+} }
Let 
$$\left(\begin{array}{cc}
			\displaystyle{L(\lambda)} & S(\lambda) \\
			T(\lambda) & V(\lambda)
		\end{array}\right) \ ,$$
be the matrix of $R(\lambda) = (\lambda I - \tilde{T})^{-1}$ 
with respect to the decomposition $\tilde{H} = H \oplus H'$. Using the 
fact that $\tilde{T}$ is an extension of $T$ and the equality 
$$\left(\begin{array}{cc}
			\displaystyle{L(\lambda)} & S(\lambda) \\
			T(\lambda) & V(\lambda)
		\end{array}\right) \left( \lambda I - \tilde{T}\right) {h \choose 0} 
= {h \choose 0} \ , h \in H \ ,$$
we obtain
$$L(\lambda)(\lambda I - T)h = h$$
for all $h \in H$. We identify $h \in H$ with ${h \choose 0} \in 
\tilde{H}$.

We have,
$$\frac{L(\lambda) - L(\mu)}{\mu - \lambda}{h \choose 0} 
= P_{H}\frac{R(\lambda) - R(\mu)}{\mu - \lambda}{h \choose 0} = 
P_{H}R(\lambda)R(\mu){h \choose 0} \ .$$
Therefore
$${\frac{L(\lambda)h - L(\mu)h}{\mu - \lambda} \choose 0} 
= \left(\begin{array}{cc}
			\displaystyle{1} & 0 \\
			0 & 0
		\end{array}\right) \left(\begin{array}{cc}
			\displaystyle{L(\lambda)} & S(\lambda) \\
			T(\lambda) & V(\lambda)
		\end{array}\right) \left(\begin{array}{cc}
			\displaystyle{L(\mu)} & S(\mu) \\
			T(\mu) & V(\mu)
		\end{array}\right) {h \choose 0}$$
for all $h \in H$, yielding
$$\frac{L(\lambda)h - L(\mu)h}{\mu - \lambda} 
= L(\lambda)L(\mu)h + S(\lambda)T(\mu)h \ ,$$
for all $h \in H$.

Since $S(\lambda)T(\lambda) = 0$, we obtain $L'(\lambda) = - 
L(\lambda)^2$, for all 
$\lambda \in \Omega$. Then $L(\lambda)$ is a left resolvent (cf. for 
instance 
\cite[Theorem 2.7]{BaMb}).
\qed   

\subsection{Proof of Theorem \ref{ext}.}

Without loss of generality one may assume that $\Omega \subseteq \D$.

Suppose Hilbert space $H'$, operators $S, U(\lambda) \in L(H')$, 
$\lambda \in \Omega$, 
$\Omega \ni \lambda \to K(\lambda) \in H'$ 
and linear functional $r : H' \to \C$ satisfying 
 \begin {enumerate}
			\item $(\lambda I - S)U(\lambda) = I$ ;   $\lambda \in \Omega$
			\item $U(\lambda)(\lambda I - S) = I - r(\cdot)K(\lambda)$ ; 
$\lambda \in \Omega$
		\item $r(K(\lambda) = 1$ ; $\lambda \in \Omega$ (normalization)
\end{enumerate} 
are fixed. This set is nonvoid as the following example shows.

\bigskip

\noindent {\bf Example.} Let $H' = \ell^2$ and consider $S$ the backward 
shift operator, 
$U(\lambda) = (\lambda U - I)^{-1}U, \; \;  \lambda \in 
\Omega$, where $U$ is the upward shift, $r(\cdot) = (\cdot,e_1)$, where 
$e_1 = (1, 0, 0, \dots)$ and 
$K(\lambda) = - (\lambda U - I)^{-1}e_1 = (1, \lambda , \lambda^2 , 
\lambda^3 , \dots )$. Then conditions 
1) -- 3) above are satisfied. 

Indeed, we have $SU = I$ and $US = I - r(\cdot)e_1$, and thus 
$$(\lambda I - S)U(\lambda) = (\lambda I - S)U(\lambda U - I)^{-1}  = I$$
and 
\begin{eqnarray*}
U(\lambda)(\lambda I - S) & = & (\lambda U - I)^{-1}U(\lambda I - S) \\
						&	= & (\lambda U - I)^{-1}\{\lambda U - I  + r(\cdot)e_1\} \\
		&	= & I + r(\cdot)(\lambda U - I)^{-1}e_1 \\
																					&	= & I - r(\cdot)K(\lambda).
\end{eqnarray*}
Condition 3) is easily verified.

\qed
\bigskip

The following properties will be useful :

\noindent {\bf Properties.}
\begin{itemize}
\item[a)] $r(u) = 0$ for all $u \in {\mbox{Im}}(U(\lambda))$.
\item[b)] $(\lambda I - S)K(\lambda) = 0$.
\end{itemize}

{\it Proof of a).\/} The condition 2) in the point $U(\lambda)h'$ gives, 
using also 
the condition 1),
$$U(\lambda)h' = U(\lambda)(\lambda I - S)U(\lambda)h' = U(\lambda)h' - 
r(U(\lambda)h')K(\lambda).$$
Therefore $r(U(\lambda)h')K(\lambda) = 0$. Using the normalization 
condition we obtain 
$r(U(\lambda)h') = 0$. \qed   

{\it Proof of b).\/} Using the second condition 2) and the normalization 
one, we get 
$$U(\lambda)(\lambda I - S)K(\lambda) = K(\lambda) - K(\lambda) = 0.$$
Since $U(\lambda)$ is left-invertible by condition 1), we obtain the 
desired conclusion. \qed

To start off the proof of Theorem \ref{ext}, we remark that $h - (\lambda 
I - T)L(\lambda) 
\in N(L(\lambda))$ for all $h \in H$. We obtain $H = Im(\lambda I - T)$ 
\.+ $N(L(\lambda))$, as 
an algebraic direct sum, for all $ \lambda \in \Omega$. Indeed, if 
$h = (\lambda I - T)u \in N(L(\lambda))$, then
$$u = L(\lambda)(\lambda I - T)u = L(\lambda)h = 0$$
and thus $h = 0$.

It is easy to see that $\dim N(L(\lambda))$ is constant for all $\lambda 
\in \Omega$. 
We suppose in this part of the proof that $\dim N(L(\lambda)) = \infty$ 
for all 
$\lambda \in \Omega$. If $\dim N(L(\lambda)) = m < \infty$, the proof will 
be essentially the same 
(in fact even simpler since we have to deal with finite sums).
Let $e_{j}(\lambda)$, $j = 1, \dots , m$, be  an orthonormal 
basis of $N(L(\lambda))$. Every $h \in H$ can be written as 
$$h = (\lambda I - T)L(\lambda)h + 
\sum_{j=1}^{\infty}c_j^{\lambda}(h)e_j(\lambda),$$
where $c_j^{\lambda}(h) \in \C$. Since $L(\lambda)(\lambda I - T) = I$, we 
get 
$$(\lambda I - T)u = (\lambda I - T)u + 
\sum_{j=1}^{\infty}c_j^{\lambda}((\lambda I - T)u)e_j(\lambda).$$
Therefore $c_j^{\lambda}(h) = 0$ if $h \in {\rm{ Im}}(\lambda I - T)$. 
Since $e_j(\lambda) \in N(L(\lambda))$, 
it follows that $c_j^{\lambda}(e_k(\lambda)) = \delta_j^k$ for $j,k \geq 
1$.

\bigskip

{\it Construction of $\tilde{T}(\lambda)$.\/} Consider $H'' = H' \oplus H' 
\oplus \dots $ the infinite direct sum of $H'$ (if 
$m < \infty$, we can consider only the finite direct sum with $m$ times 
$H'$). 
Set $\tilde{H} = H \oplus H''$. Let $\tilde{T}(\lambda) \in 
L(\tilde{H})$ be defined by 
$$\tilde{T}(\lambda)(h,h'_{1},  h'_2, \dots ) = (Th - 
\sum_{j=1}^{\infty}r(h'_{j})e_{j}(\lambda), Sh'_{1}, Sh'_2, \dots ).$$
Define also $\tilde{R}(\lambda) \in L(\tilde{H})$ by  
$$\tilde{R}(\lambda)(h, h'_{1}, h'_{2}, \dots ) = 
(L(\lambda)h, c_{1}^{\lambda}(h)K(\lambda) +U(\lambda)h'_{1},
c_{2}^{\lambda}(h)K(\lambda) + U(\lambda)h'_{2}, \dots ).$$

\bigskip

{\it Proof of 1.\/} It can be easily seen that  $\tilde{T}(\lambda)$ 
and $\tilde{R}(\lambda)$ are well defined. It is also 
easy to see that $\tilde{T}(\lambda)(h,0,0, \dots ) = (Th, 0,0, \dots )$ 
for all 
$h \in H$. Therefore $\tilde{T}(\lambda)$ is an extension of $T$. 

\bigskip

{\it Proof of 2.\/} We show now that $\tilde{R}(\lambda) = (\lambda I - 
\tilde{T}(\lambda))^{-1}$. Indeed, we have
$$\tilde{R}(\lambda)(\lambda I - \tilde{T}(\lambda))(h,h'_1, \dots)$$
\begin{center}
$= \tilde{R}(\lambda)((\lambda I - T)h + 
\sum_{j=1}^{\infty}r(h'_j)e_j(\lambda), (\lambda I - S)h'_1, \dots )$
\end{center}
$$= (L(\lambda)(\lambda I - T)h,r(h'_1)K(\lambda) + U(\lambda)(\lambda I - 
S)h'_1, \dots )$$
\begin{center}
(since $L(\lambda)e_j(\lambda) = 0$ and $c_j^{\lambda}(u) = 0$ 
for all $u \in {\rm{ Im}}(\lambda I - T))$
\end{center}
$$= (h,h'_1, \dots ).$$

On the other hand,
$$(\lambda I - \tilde{T}(\lambda))\tilde{R}(\lambda)(h,h'_1, \dots )$$
\begin{center}
$= (\lambda I - \tilde{T}(\lambda))(L(\lambda)h,c_1^{\lambda}(h)K(\lambda) 
+ U(\lambda)h'_1, \dots )$
\end{center}
$$= ((\lambda I - T)L(\lambda)h + 
\sum_{j=1}^{\infty}c_j^{\lambda}(h)e_j(\lambda),(\lambda I - 
S)(c_1^{\lambda}(h)K(\lambda) 
+ U(\lambda)h'_1),\dots )$$
\begin{center}
(using Property a)
\end{center}
$$ = (h,(\lambda I - S)U(\lambda)h'_1, \dots )$$
\begin{center}
(using Property b)
\end{center}
$$= (h,h'_1, \dots ).$$

Thus $\tilde{R}(\lambda) = (\lambda I - \tilde{T}(\lambda))^{-1}$ and 
Property 2 in the Theorem is 
proved. 

\bigskip

{\it Proof of 3.\/} It is clear from the definition that, with respect to 
the decomposition $\tilde{H} = H \oplus H''$, one has
$$(\lambda I - \tilde{T}(\lambda))^{-1} = \tilde{R}(\lambda) = 
\left(\begin{array}{cc}
			\displaystyle{L(\lambda)} & 0 \\
			\ast & \ast
		\end{array}\right) \ .$$

\bigskip

{\it Proof of 4.\/} The fourth property ($L(\lambda)$ is a left resolvent 
iff $\tilde{T}(\lambda)$ is 
independent of $\lambda$) follows from the definition of 
$\tilde{T}(\lambda)$ and the 
fact that a left-inverse function $L(\lambda)$ is a left resolvent if and 
only if the kernel $N(L(\lambda))$ is independent of $\lambda$ (cf. 
\cite[Lemma 9.6]{HerII}).

\bigskip

{\it Proof of 5.\/} Suppose, for the last property, that $m < \infty$. 
According to a result of  
\cite{Sap1} and \cite[page 187]{Bar} , the orthonormal basis 
$e_j(\lambda), j = 1, \dots, m$, 
can be chosen 
such that $\lambda \to e_j(\lambda)$ is analytic. This implies that the 
map $\lambda \to \tilde{T}(\lambda)$ is 
analytic.
\qed

\end{document}